\documentclass{article}
\usepackage[english]{babel}
\usepackage{geometry,amsmath,graphicx,bbm}
\newcommand{\Zeta}{\mathrm{Z}}
\newcommand{\mathd}{\mathrm{d}}
\newcommand{\tmem}[1]{{\em #1\/}}
\newcommand{\tmop}[1]{\ensuremath{\operatorname{#1}}}
\newcommand{\tmtextit}[1]{{\itshape{#1}}}

\begin{document}

\title{The Laplacian of The Integral Of The Logarithmic Derivative of the
Riemann-Siegel-Hardy Z-function}

\author{{\hspace{3em}}Stephen Crowley}

\maketitle

\begin{abstract}
  The integral $R (t) = \pi^{- 1} \int \frac{\mathd}{\mathd t} \tmop{lnZ} (t)
  \mathd t = \pi^{- 1}  \left( \ln \zeta \left( \frac{1}{2} + i t \right) + i
  \vartheta (t) \right)$ of the logarithmic derivative of the Hardy Z function
  $Z (t) = e^{i \vartheta (t)} \zeta \left( \frac{1}{2} + i t \right)$, where
  $\vartheta (t)$ is the Riemann-Siegel theta function, and $\zeta (t)$ is the
  Riemann zeta function, is used as a basis for the construction of a pair of
  transcendental entire functions $\nu (t) = - \nu (1 - t) = - \left( \Delta R
  \left( \frac{i}{2} - i t \right) \right)^{- 1} = - G \left( \frac{i}{2} - i
  t \right)$ where $G = - (\Delta R (t))^{- 1}$ is the derivative of the
  additive inverse of the reciprocal of the Laplacian $\Delta f (t) = \ddot{f}
  (t)$ of $R (t) $ and $\chi (t) = - \chi (1 - t) = \dot{\nu} (t) = - i H
  \left( \frac{i}{2} - i t \right)$ where $H (t)$=$\dot{G} (t)$ has roots at
  the local minima and maxima of $G (t)$. When $H (t) = 0$ and $\dot{H} (t) =
  \ddot{G} (t) = \Delta G (t) > 0$, the point $t$ marks a minimum of $G (t)$
  where it coincides with a Riemann zero, i.e., $\zeta \left( \frac{1}{2} + i
  t \right) = 0$, otherwise when $H (t)^{} = 0$ and $\dot{H} (t) = \Delta G
  (t) < 0$, the point $t$ marks a local maximum $G (t)$, marking midway points
  between consecutive minima. Considered as a sequence of distributions or
  wave functions, $\nu_n (t) = \nu (1 + 2 n + 2 t)$ converges to $\nu_{\infty}
  (t) = \lim_{n \rightarrow \infty} \nu_n (t) = \sin^2 (\pi t)$ and $\chi_n
  (t) = \chi (1 + 2 n + 2 t)$ to \ $\chi_{\infty} (t) = \lim_{n \rightarrow
  \infty} \chi_n (t) = - 8 \cos (\pi t) \sin (\pi t)$ 
\end{abstract}

{\tableofcontents}

\section{Derivations}

\subsection{Standard Definitions}

Let $\zeta (t)$ be the Riemann zeta function
\begin{equation}
  \begin{array}{lll}
    \zeta (t) & = \sum_{n = 1}^{\infty} n^{- s} & \forall \tmop{Re} (s) > 1\\
    & = (1 - 2^{1 - s}) \sum_{n = 1}^{\infty} n^{- s} (- 1)^{n - 1} & \forall
    \tmop{Re} (s) > 0
  \end{array}
\end{equation}
and $\vartheta (t)$ be Riemann-Siegel vartheta function $\vartheta (t)$

\begin{equation}
  \vartheta (t) = - \frac{i}{2} \left( \ln \Gamma \left( \frac{1}{4} + \frac{i
  t}{2} \right) - \ln \Gamma \left( \frac{1}{4} - \frac{i t}{2} \right)
  \right) - \frac{\ln (\pi) t}{2}
\end{equation}
where $\arg (z) = \frac{\ln (z) - \ln (\bar{z})}{2 i}$ and $\overline{\Gamma
(z)} = \Gamma (\bar{z})$. The Hardy $Z$ function{\cite{HardyZ}} can then be
written as
\begin{equation}
  Z (t) = e^{i \vartheta (t)} \zeta \left( \frac{1}{2} + i t \right) \label{Z}
\end{equation}
which can be mapped isometrically back to the $\zeta$ function
\begin{equation}
  \zeta (t) = e^{- i \vartheta \left( \frac{i}{2} - i t \right)} Z \left(
  \frac{i}{2} - i t \right) \label{Zz}
\end{equation}
due to the isometry
\begin{equation}
  t = \frac{1}{2}_{} + i \left( \frac{i}{2} - i t \right)
\end{equation}
of the Mobius transforms\footnote{Thanks to Matti Pitk{\"a}nen for pointing
out this is a Mobius transform pair, among other things} $f (t) = \frac{a t +
b}{c t + d} = \left(\begin{array}{cc}
  a & b\\
  c & d
\end{array}\right)$ with
\begin{equation}
  \left(\begin{array}{cc}
    - i & \frac{i}{2}\\
    0 & 1
  \end{array}\right) \tmop{and} \tmop{its} \tmop{inverse}
  \left(\begin{array}{cc}
    i & \frac{1}{2}\\
    0 & 1
  \end{array}\right)
\end{equation}
making possible the Riemann-Siegel-Hardy correspondence. Furthermore, let $S
(t)$ be argument of $\zeta$ normalized by $\pi$ defined by
\begin{equation}
  \begin{array}{ll}
    S (t) & = \pi^{- 1} \arg \left( \zeta \left( \frac{1}{2} + i t \right)
    \right)\\
    & = - \frac{i}{2 \pi} \left( \ln \zeta \left( \frac{1}{2} + i t \right) -
    \ln \zeta \left( \frac{1}{2} - i t \right) \right)
  \end{array} \label{S}
\end{equation}
The B{\"a}cklund counting formula gives the exact number of zeros on the
critical strip up to level $t$, not just on the criticial line,
\begin{equation}
  N (t) = \tmop{Im} (R (t)) = \frac{\vartheta (t)}{\pi} + 1 + S (t) \label{N}
\end{equation}
The relationship between the functions $N (t)$, $S (t)$, and $Z (t)$ is
demonstrated by
\begin{equation}
  \ln \zeta \left( \frac{1}{2}_{} + i t \right) = \ln | Z (t) | + i \pi S (t)
\end{equation}
These formulas are true independent of the Riemann hypothesis which posits
that all complex zeros of $\zeta (s + i t)$ have real part $s = \frac{1}{2}$.
{\cite[Corrollary 1.8 p.13]{HardyZ}}

\subsection{The Logarithmic Derivative of $Z (t)$ and its Integral}

Let $Q (t)$ be the logarithmic derivative of $Z (t)$ given by
\begin{equation}
  \begin{array}{ll}
    Q (t) & = \frac{\dot{\Zeta} (t)}{Z (t)}\\
    & = \frac{\mathd}{\mathd t} \tmop{lnZ} (t)\\
    & = i \left( \frac{\dot{\zeta} \left( \frac{1}{2} + i t \right)}{\zeta
    \left( \frac{1}{2} + i t \right)} + \frac{\Psi \left( \frac{1}{4} -
    \frac{i}{2 t} \right) + \Psi \left( \frac{1}{4} + \frac{i}{2 t} \right) -
    2 \ln (\pi)}{4} \right)\\
    & = i \left( \frac{\dot{z} (t)}{z (t)} + \frac{g^- (t) + g^+ (t) - 2 \ln
    (\pi)}{4} \right)
  \end{array} \label{Q}
\end{equation}
where
\begin{equation}
  \Psi (x) = \frac{\mathd}{\mathd x} \ln (\Gamma (x)) =
  \frac{\frac{\mathd}{\mathd x} \Gamma (x)}{\Gamma (x)}
\end{equation}
is the digamma function, the logarithmic derivative of the $\Gamma$ function
and $z (t)$ and $g^{\pm} (t)$ have been introduced to simplify the
expressions. The function $Q (t)$ has singularities at $\pm \frac{i}{2} (4 n -
3)$ with residues
\begin{equation}
  \underset{t = \frac{i}{2} (4 n - 3)}{\tmop{Res} (Q (t))} = \frac{- 8 \zeta
  (2 - 2 n) n^2 + 4 \zeta (2 - 2 n) n}{16 \zeta (2 - 2 n) n^2 - 8 \zeta (2 - 2
  n) n} = - \frac{1}{2}
\end{equation}
and
\begin{equation}
  \underset{t = - \frac{i}{2} (4 n - 3)}{\tmop{Res} (Q (t))} = \frac{- 2 \zeta
  (2 n - 1) + 4 \zeta (2 n - 1) n}{8 \zeta (2 n - 1) n^{} - 4 \zeta (2 n - 1)}
  = \frac{1}{2}
\end{equation}
Now, the integral of the logarithmic derivative of Z is defined by
\begin{equation}
  \begin{array}{cl}
    R (t) & = \pi^{- 1} \int Q (t) \mathd t\\
    & = \pi^{- 1} \int \frac{\dot{\Zeta} (t)}{Z (t)} \mathd t\\
    & = \pi^{- 1} \int \frac{\mathd}{\mathd t} \tmop{lnZ} (t) \mathd t\\
    & = \pi^{- 1} \int i \left( \frac{\dot{z} (t)}{z (t)} + \frac{g^- (t) +
    g^+ (t) - 2 \ln (\pi)}{4} \right) \mathd t\\
    & = \pi^{- 1} \int i \left( \frac{\dot{\zeta} \left( \frac{1}{2} +
    \tmop{it} \right)}{\zeta \left( \frac{1}{2} + i t \right)} + \frac{\Psi
    \left( \frac{1}{4} - \frac{i}{2 t} \right) + \Psi \left( \frac{1}{4} +
    \frac{i}{2 t} \right) - 2 \ln (\pi)}{4} \right) \mathd t\\
    & = \pi^{- 1}  \left( \ln \zeta \left( \frac{1}{2} + i t \right) + i
    \vartheta (t) \right)\\
    & = \pi^{- 1}  \left( \ln \zeta \left( \frac{1}{2} + i t \right) +
    \frac{i}{2} \left( i \left( \ln \Gamma \left( \frac{1}{4} - \frac{i t}{2}
    \right) - \ln \Gamma \left( \frac{1}{4} + \frac{i t}{2} \right) \right) -
    \ln (\pi) t \right) \right)
  \end{array}
\end{equation}

\subsection{The Laplacian}

The Laplacian $= \nabla \cdot \nabla = \nabla^2$ is a differential operator
which corresponds to the {\tmem{divergence}} of the {\tmem{gradient}} of a
function $f (x)$ on a Euclidean space $x \in X \underline{\subset}
\mathbbm{R}^d$ and is denoted by
\begin{equation}
  f (x) = \sum_{i = 1}^d \frac{\partial}{\partial x^2_i} f (x_i) \forall x \in
  X \subseteq \mathbbm{R}^d
\end{equation}
and is simply the second derivative of a function when $d = 1$
\begin{equation}
  f (x) = \ddot{f} (x) = \frac{\mathd}{\mathd x^2} f (x) \forall x \in X
  \subseteq \mathbbm{R}
\end{equation}
Let $G (t)$ be the additive inverse of the reciprocal(also known as
multiplicative inverse) of the Laplacian of $R (t)$ interpreted as a partition
function
\begin{equation}
  \begin{array}{ll}
    G (t) & = - \frac{1}{\Delta R (t)}\\
    & = \frac{8 z (t)^2 \pi}{z (t)^2 (\dot{g}^- (t) - \dot{g}^{+_{}} (t)) - 8
    (z (t) \ddot{z} (t) + \dot{z} (t)^2)}
  \end{array} \label{Gpart}
\end{equation}

\begin{figure}[h]
  \begin{tabular}{l}
    \resizebox{6in}{2in}{\includegraphics{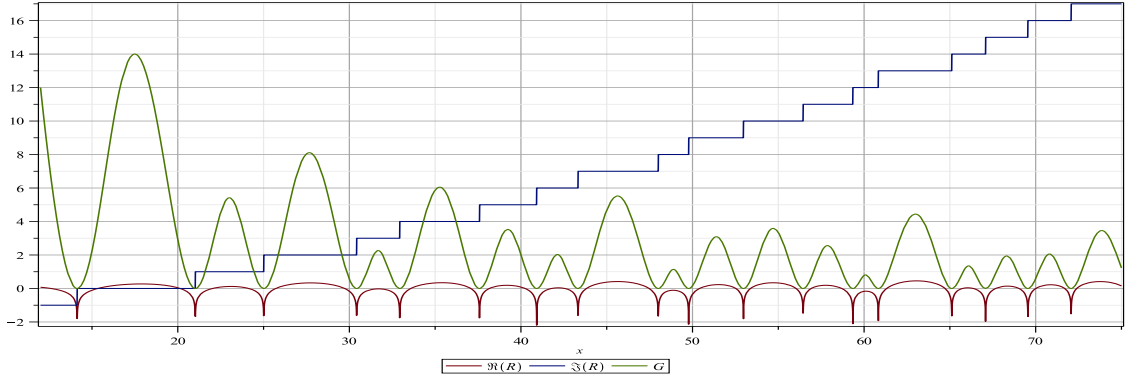}}
  \end{tabular}
  \caption{The Real and Imaginary Parts of $R (t)$ compared with $G (t)$}
\end{figure}

Now, let $H (t)$ be the derivative of $G (t)$ given by
\begin{equation}
  \begin{array}{ll}
    H (t) & = \dot{G} (t)\\
    & = \frac{\mathd}{\mathd t} \left( - \frac{1}{\Delta R (t)} \right)\\
    & = - 4 i \pi z (t) \frac{z (t)^3 \Delta g (t) + 8 (2 z (t)^2  \dddot{z}
    (t) - 6 z (t) \dot{z} (t)  \ddot{z} (t) - 4 \dot{z} (t)^3)}{(z (t)^2 + 8
    (z (t) \ddot{z} (t) - \dot{z}  (t)^2))^2}
  \end{array}
\end{equation}
When $H (t)^{} = 0$ and $\dot{H} (t) = \ddot{G} (t) = \Delta G (t) > 0$, the
point $t$ marks a minimum of $G (t)$ where it coincides with a Riemann zero,
i.e., $\zeta \left( \frac{1}{2} + i t \right) = 0$, otherwise when $H (t)^{} =
0$ and $\dot{H} (t) = \Delta G (t) < 0$, the point $t$ marks a local maximum
$G (t)$, marking midway points between consecutive minima.

\begin{figure}[h]
  \resizebox{6in}{1.5in}{\includegraphics{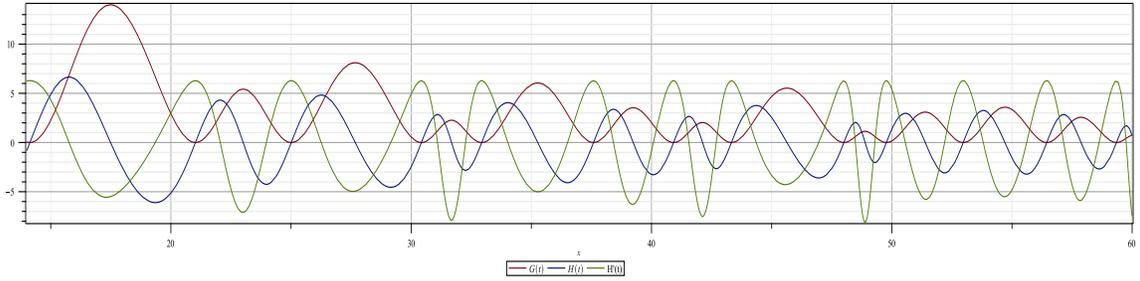}}
  \caption{$G (t), H (t) \tmop{and} \dot{H} (t)$}
\end{figure}

With the identity $t = \frac{1}{2}_{} + i \left( \frac{i}{2} - i t \right)$,
define $\nu (t)$ as the Mobius transform of the partition function
\begin{equation}
  \nu (t) = - \left( \Delta R \left( \frac{i}{2} - i t \right) \right)^{- 1} =
  - G \left( \frac{i}{2} - i t \right)
\end{equation}
then $\nu (t)$ has zeros at the positive odd integers, zero, and negative even
integers. In the same way, let
\begin{equation}
  \begin{array}{ll}
    \chi (t) & = \dot{\nu} (t)\\
    & = - i H \left( \frac{i}{2} - i t \right)\\
    & = - 4 i \pi \zeta (t)  \frac{\zeta (t)^3 \left( \ddot{\Psi} \left(
    \frac{t}{2} \right) - \ddot{\Psi} \left( \frac{1}{2} - \frac{t}{2} \right)
    \right) - 8 (6 \zeta (t) \dot{\zeta} (t) \ddot{\zeta} (t) + 2 \zeta
    (t)^{^2} \dddot{\zeta} (t) + 4 \dot{\zeta} (t)^3)}{\zeta (t)^2 \left(
    \left( \dot{\Psi} \left( \frac{t}{2} \right) - \dot{\Psi} \left(
    \frac{1}{2} - \frac{t}{2} \right) \right) + 8 (\zeta (t) \ddot{\zeta} (t)
    - \dot{\zeta} (t)^2) \right)^2}
  \end{array}
\end{equation}
which also has zeros on the real line at the positive odd integers, zero, and
the negative even integers.
\begin{equation}
  \begin{array}{rll}
    \lim_{t \rightarrow 2 n - 1} \nu (t) & = \lim_{t \rightarrow 2 n - 1} \chi
    (t) & = 0 \forall n > 0\\
    \lim_{t \rightarrow 0} \nu (t) & = \lim_{t \rightarrow 0} \chi (t) & = 0\\
    \lim_{t \rightarrow - 2 n} \nu (t) & = \lim_{t \rightarrow - 2 n} \chi (t)
    & = 0 \forall n < 0
  \end{array}
\end{equation}
Both $\nu (t)$ and $\chi (t)$ satisify similiar functional equations
\begin{equation}
  \nu (t) = \nu (1 - t)
\end{equation}
and
\begin{equation}
  \chi (t) = - \chi (1 - t)
\end{equation}
So the even $\mu (t)$ and odd $\psi (t)$ transcendental entire functions can
be defined
\begin{equation}
  \mu (t) = \nu \left( t + \frac{1}{2} \right)
\end{equation}
\begin{equation}
  \psi (t) = \chi \left( t + \frac{1}{2} \right)
\end{equation}
which satisify the functional symmetries
\begin{equation}
  \mu (t) = \mu (- t)
\end{equation}
and
\begin{equation}
  \psi (t) = - \psi (- t)
\end{equation}
The function $\chi (t)$ has as a subset of its roots, the roots of the Rieman
zeta function $\zeta (t)$, the converse is not true, since $\chi (t)$ is a
function of $\zeta (t)$ and its first, second, and third derivatives.
\begin{equation}
  \{ t : \zeta (t) = 0 \} \subset \{ t : \chi (t) = 0 \}
\end{equation}
Let
\begin{equation}
  \chi_n (t) = \chi (1 + 2 n + 2 t)
\end{equation}
and
\begin{equation}
  \nu_n (t) = \nu (1 + 2 n + 2 t)
\end{equation}
which both satisify the symmetries
\begin{equation}
  \chi_{- \frac{1}{2}} \left( \frac{t}{2} \right) = \chi (t)
\end{equation}
\begin{equation}
  \nu_{- \frac{1}{2}} \left( \frac{t}{2} \right) = \nu (t)
\end{equation}
as well as
\begin{equation}
  \chi_n \left( - \frac{1}{2} \right) = \chi (2 n)
\end{equation}
\begin{equation}
  \nu_{_n} \left( - \frac{1}{2} \right) = \nu (2 n)
\end{equation}
The sequence of wave functions $\nu_n (t)$ converges, thanks to quantum
ergodicity, to
\begin{equation}
  \nu_{\infty} (t) = \lim_{n \rightarrow \infty} \nu_n (t) = \sin^2 (\pi t)
\end{equation}
which has a limiting maximum

\begin{equation}
  \begin{array}{ll}
    \lim_{n \rightarrow \infty} \max_{0 < t < 1} \nu_n (t) & = \max_{0 < t <
    1} \nu_{\infty} (t)\\
    & = \lim_{n \rightarrow \infty} \nu (2 n)\\
    & = \frac{8}{\pi}\\
    & \cong 2.546479089470 \ldots
  \end{array}
\end{equation}
and the associated differential is
\begin{equation}
  \begin{array}{cl}
    \chi_{\infty} (t) & = \lim_{n \rightarrow \infty} \chi_n (t)\\
    & = \frac{4}{\pi} \frac{\mathd}{\mathd t} \nu_{\infty} (t)\\
    & = \frac{4}{\pi} \frac{\mathd}{\mathd t} \sin^2 (\pi t)\\
    & = - 8 \sin (\pi t) \cos (\pi t)
  \end{array}
\end{equation}
It is worth mentioning that the pair-correlation function for the zeros of
$\zeta (s)$ is
\[ r_2 (x) = 1 - \frac{\sin^2 (\pi x)^2}{\pi^2 x^2} \]
if the Riemann hypothesis is true. {\cite{prz}} The point $n = 1$ is where
$\chi_n \left( - \frac{1}{2} \right)$ and $\nu_n \left( - \frac{1}{2} \right)$
attain their greatest values among the integers.

\begin{table}[h]
  \begin{tabular}{l}
    $\begin{array}{ccc}
      &  & 
    \end{array}$
  \end{tabular}\begin{tabular}{|c|c|c|}
    \hline
    $n$ & $\nu_n \left( - \frac{1}{2} \right)$ & $\chi_n \left( - \frac{1}{2}
    \right)$\\
    \hline
    $1$ & 117.43532857805377782 & 9447.7593604718560\\
    \hline
    $2$ & 003.03320654562255410 & 0000.2816402783351\\
    \hline
    $3$ & 002.77176105375846239 & 0000.0589526080385\\
    \hline
    $4$ & 002.69653220844944185 & 0000.0240373109008\\
    \hline
    5 & 002.66095375057354766 & 0000.0132398305603\\
    \hline
    $6$ & 002.63970550787458574 & 0000.0088555925215\\
    \hline
    $7$ & 002.62532107269483772 & 0000.0060558883634\\
    \hline
    {\ldots} & $\ldots \ldots \ldots \ldots \ldots \ldots$ & $\ldots \ldots
    \ldots \ldots \ldots \ldots \ldots$\\
    \hline
    $\infty$ & $\frac{8}{\pi}$ & 0\\
    \hline
  \end{tabular}
  \caption{Even values of $\nu (n)$ and $\chi (n)$ }
\end{table}

This phenomena of having the first wave function having a much larger size
than all of the remainders is mentioned in {\cite[Theorem 22]{knauf99}} The
convergence of $\lim_{n \rightarrow \infty} \nu_n (t) \rightarrow \nu_{\infty}
(t)$ and $\lim_{n \rightarrow \infty} \chi_n (t) \rightarrow \chi_{\infty}
(t)$ might be interpreted as a manifestation of quantum ergodicity {\cite[C
p.16]{prz}}. There is an essential singularity{\cite[54 p.169]{cva}} at $s_0
\cong 1.98757823 \ldots$
\begin{equation}
  \lim_{s \rightarrow s_0} \chi (s) = \infty
\end{equation}
We also have the limits
\begin{equation}
  \lim_{s \rightarrow \{ 0, 1 \}} \dot{\chi} (s) = 4 \pi
\end{equation}
The integral over this sine-squared kernel is
\begin{equation}
  \int_0^1 \sin^2 (\pi t) \mathd t = \frac{1}{2}
\end{equation}
whereas
\begin{equation}
  \frac{}{} \int_0^1 \frac{\nu (t)}{\nu \left( \frac{1}{2} \right)} \mathd t
  \cong .46693755153559653755 \ldots
\end{equation}

\begin{figure}[h]
  \resizebox{2in}{1.75in}{\includegraphics{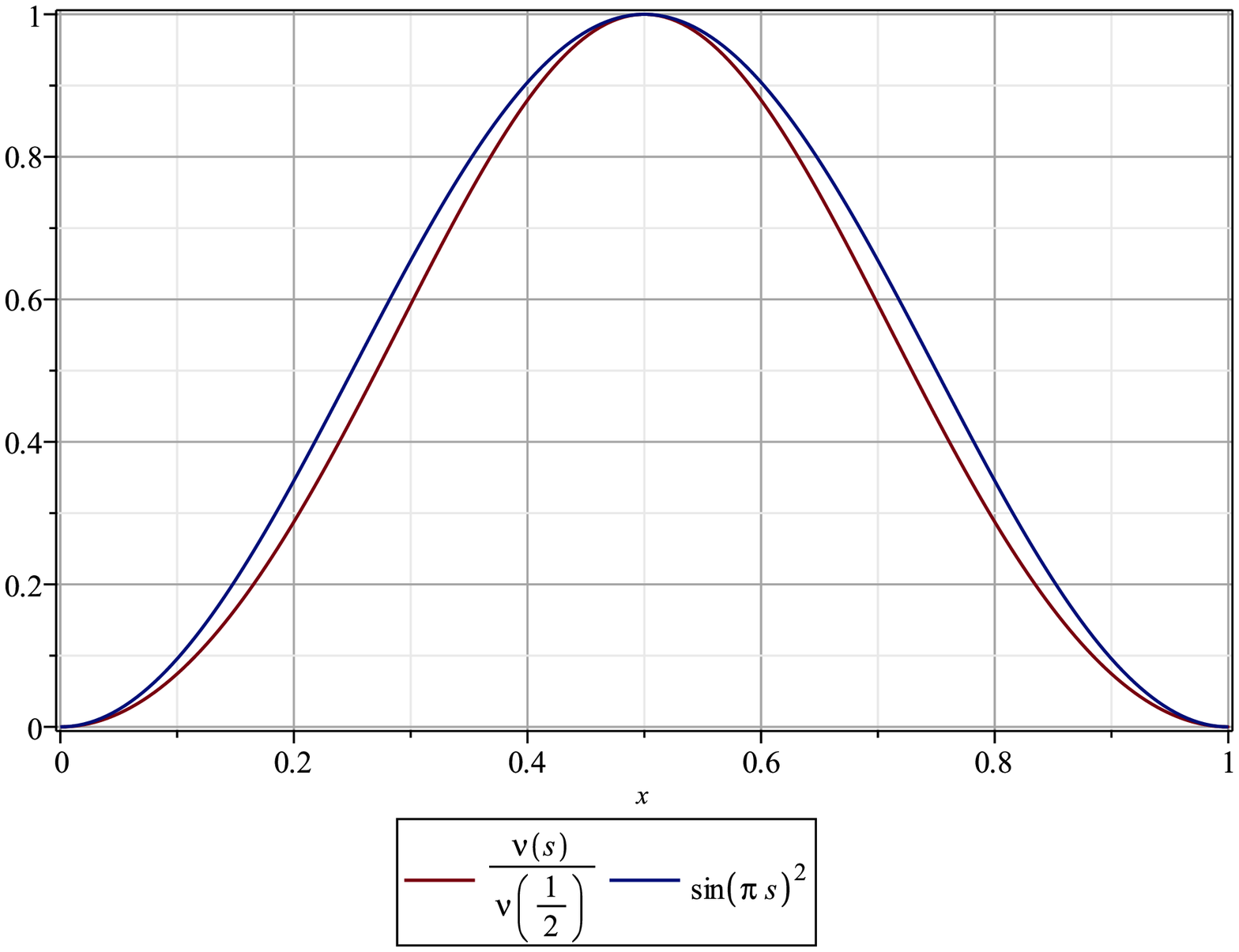}}\resizebox{2.5in}{1.75in}{\includegraphics{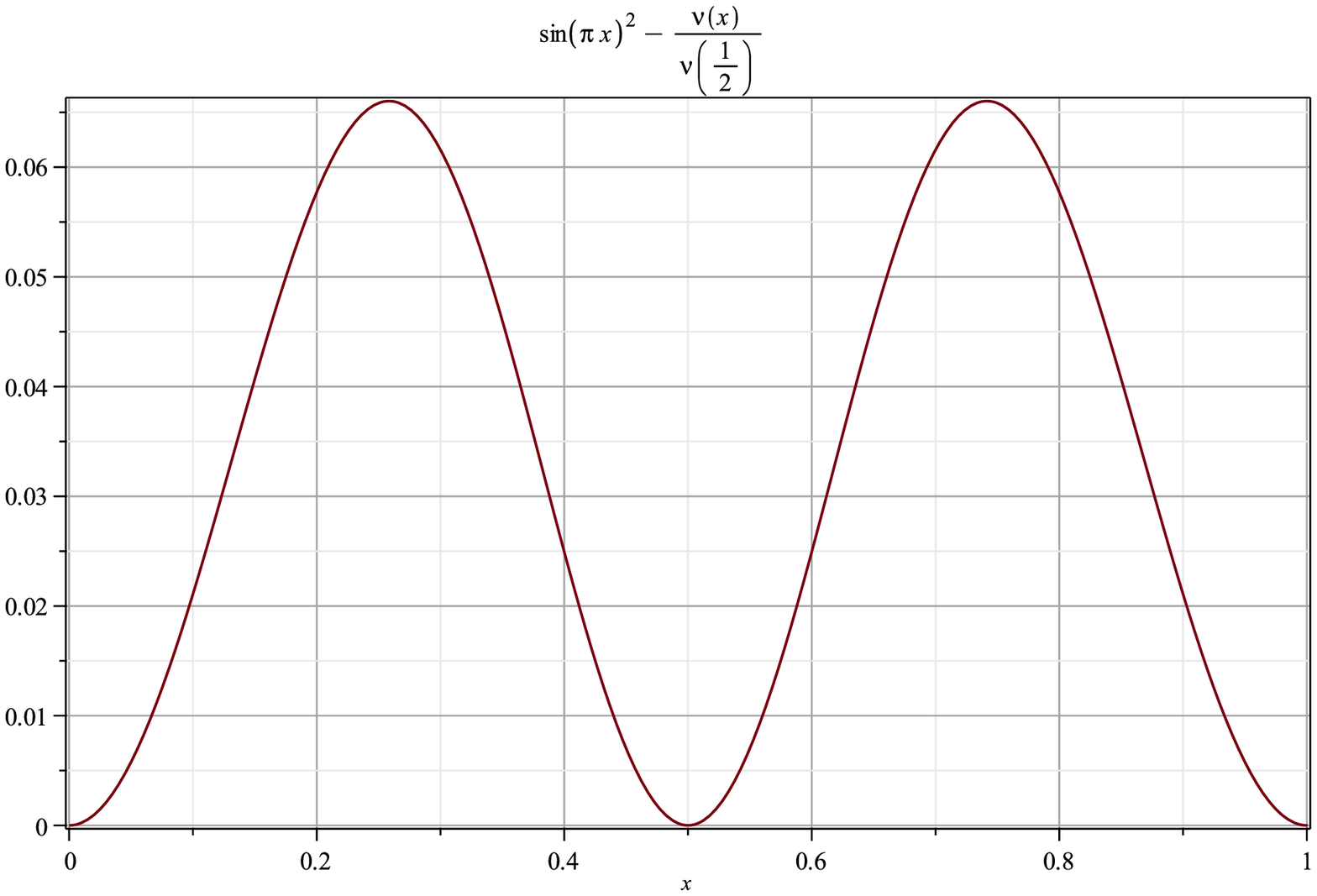}}
  
  \resizebox{2.5in}{1.75in}{\includegraphics{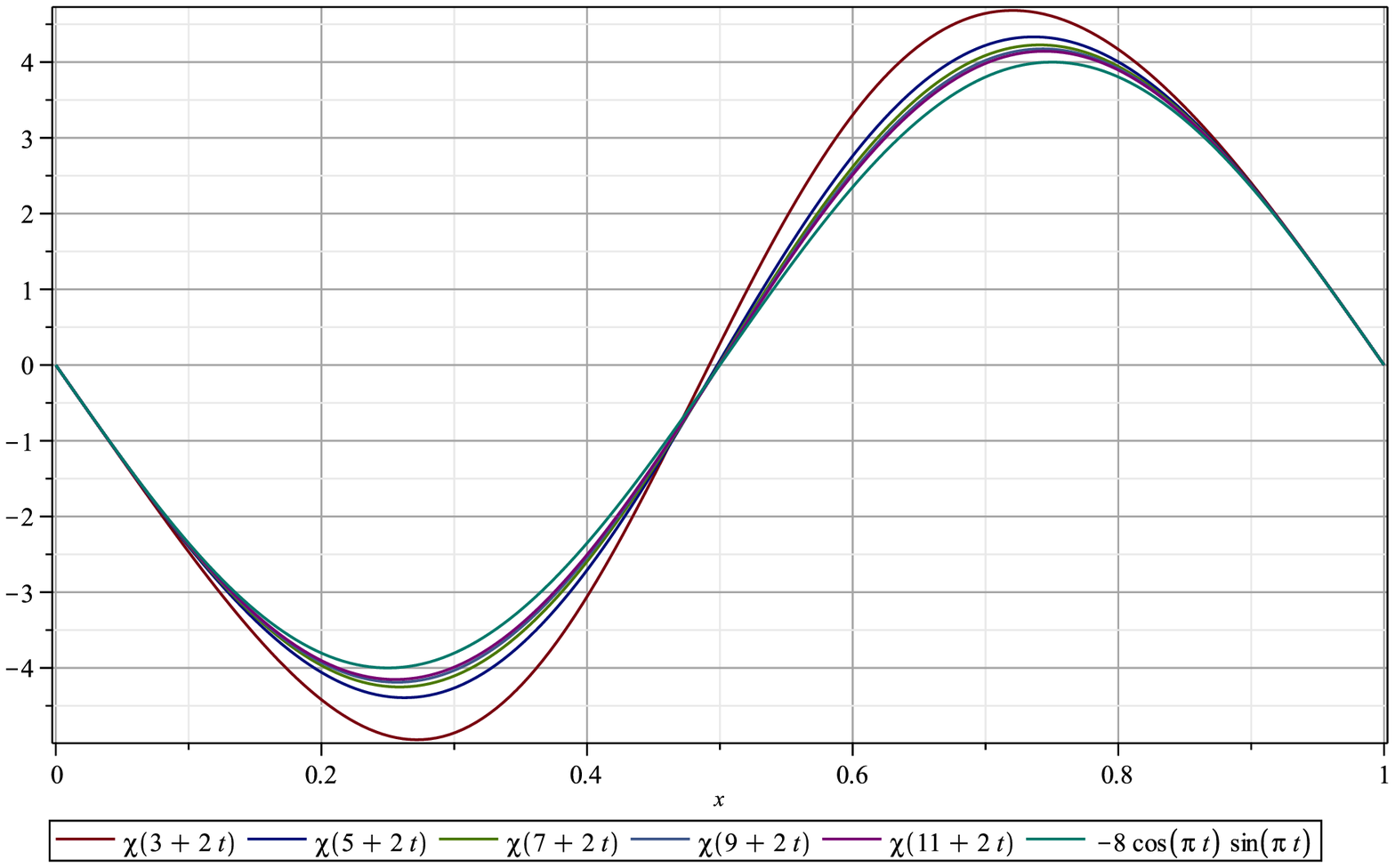}}\resizebox{2.5in}{1.75in}{\includegraphics{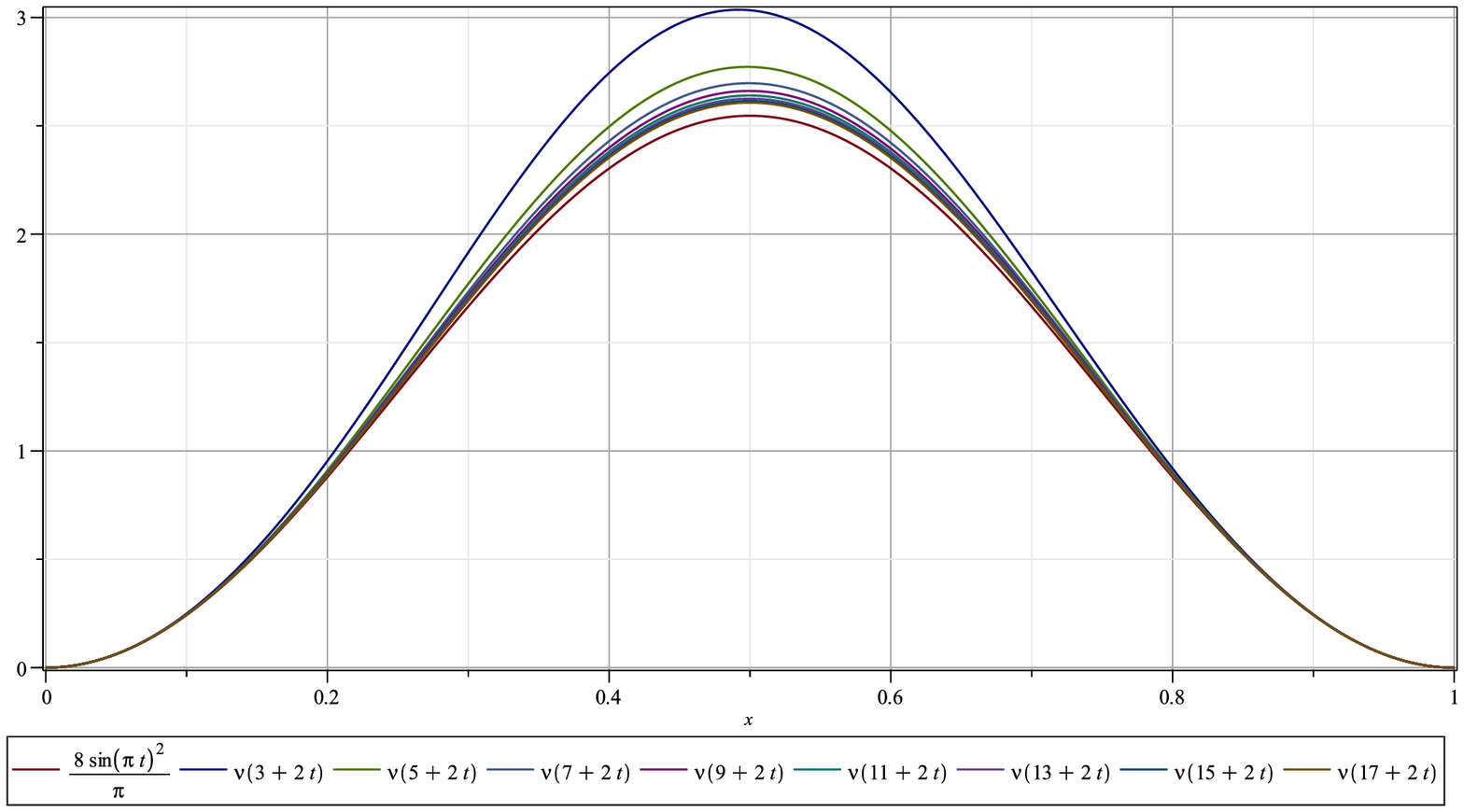}}
  \caption{{\tmem{Top Left}}: The function $\nu (s)$ normalized by its maximum
  value on $[0, 1]$ {\tmem{Top Right}}: The function $\nu (s)$ normalized by
  its maximum value and subtracted from $\sin^2 (\pi x)$ {\tmem{Bottom Left}}:
  Convegence of $_{} \chi_n (t) \rightarrow - 8 \cos (\pi t) \sin^{} (\pi t)$
  \ {\tmem{Bottom Right}}: Convegence of $_{} \nu_n (t) \rightarrow \sin^2
  (\pi t) \frac{\pi}{8}$}
\end{figure}

\section{Appendix}

\subsection{Wave Mechanics}

\subsection{The One-Dimensional Wave Equation}

The one-dimensional wave equation is
\begin{equation}
  - \frac{\hbar^2}{2 m} \frac{\partial^2 }{x^2} u (x, t) + V (x, t) u (t, x) =
  - \frac{\hbar}{i} \frac{\partial}{\partial t} u (x, t)
\end{equation}
{\cite[II.A.4 p.25]{pqm}}{\cite[4.2.1 p.]{MonsterMoonshine}}{\cite[Ch. 8
p.147]{Psm}}

\subsection{The Poisson Bracket and Lagrangian Mechanics}

``The'' Poisson bracket, expressed with Einstein summation convention (for the
repeated index $i$)
\begin{equation}
  \{ A, B \} \equiv \frac{\partial A}{\partial p_i} \frac{\partial B}{\partial
  q_i} - \frac{\partial B}{\partial p_i} \frac{\partial A}{\partial q_i}
\end{equation}
has the {\tmem{antisymmetry}} property
\begin{equation}
  \{ A, B \} = - \{ B, A \}
\end{equation}
and the so-called {\tmem{Jacobi identity}}
\begin{equation}
  \{ A, \{ B, C \} \} - \{ B, \{ A, C \} \} + \{ C, \{ A, B \} \} = 0
\end{equation}
Two quantitities $A$ and $B$ are said to {\tmem{commute}} if their
$\tmop{Poisson} \tmop{bracket}$ $\{ A, B \}$ vanishes, that is, $\{ A, B \} =
0$. Hamilton's equations of motion for the system
\begin{equation}
  \dot{p}_i = - \frac{H}{q_i}
\end{equation}
\begin{equation}
  \dot{q}_i = \frac{\partial H}{\partial q_i}
\end{equation}
where $H$ is a Legendre-transformed function of the Lagrangian called the
Hamiltonian
\begin{equation}
  H \equiv \frac{\partial L}{\partial \dot{q}_i} \dot{q}_i - L (q_i,
  \dot{q}_i, t)
\end{equation}
whose value for any given time $t$ gives the {\tmem{energy}}
\begin{equation}
  \mathcal{A} [q_i] = \int_{t_a}^{t_b} L (q_i, \dot{q}_i, t) \mathd t
\end{equation}
of the system where $L (q_i, \dot{q}_i, t)$ is the Lagrangian of the system
and
\begin{equation}
  q_i (t) = q_i^{\tmop{cl}} (t) + \delta q_i (t)
\end{equation}
is an arbitrary path where $q_i^{\tmop{cl}} (t)$ is the {\tmem{classical
orbit}} or {\tmem{classical path}} of the system and
\begin{equation}
  \delta q_i (t_a) \equiv q_i (t) - q^{\tmop{cl}}_i (t)
\end{equation}
{\cite[1.1]{kleinert2004path}}

\subsubsection{The Euler-Lagrange equation}

The Euler-Lagrange equation
\begin{equation}
  \frac{\mathd}{\mathd t} \left( \frac{\partial L}{\partial \dot{q}_i} \right)
  = \frac{\partial L}{\partial q_i}
\end{equation}
indicates that the action $S$ given by
\begin{equation}
  S = \int_{t_1}^{t_2} L (t) \mathd t
\end{equation}
where
\begin{equation}
  L (t) = T - V (x (t))
\end{equation}
is the {\tmem{Lagrangian}} and
\begin{equation}
  T (t) = \frac{1}{2} m \dot{x} (t)^2
\end{equation}
is the {\tmem{kinetic energy}} \ which is stationary for the physical
solutions $q_i (t)$.{\cite[4.2.1]{MonsterMoonshine}}

\subsubsection{Quantum Mechanics of General Lagrangian Systems}

The coordinate transformation
\begin{equation}
  x^i = x^i (q^{\mu})
\end{equation}
implies the relation
\begin{equation}
  \partial_{\mu} = \frac{\partial}{\partial q^{\mu}} = e^i_{\mu} (q)
  \partial_i \label{dm}
\end{equation}
between the derivatives $\partial_{\mu}$ and
\begin{equation}
  \partial_i \equiv \frac{\partial}{\partial x^i}
\end{equation}
where
\begin{equation}
  e^i_{\mu} (q) \equiv \partial_{\mu} x^i (q)
\end{equation}
is a transformation matrix called the {\tmem{basis p-ad}} where $p$ is the
prefix corresponding to $n$, the dimension of $x$, monad when $n = 1$, dyad
when $n = 2$, triad when $n =$3, and so on. Let
\begin{equation}
  e_i^{\mu} (q) = \frac{\partial q^{\mu}}{\partial x^i}
\end{equation}
be the inverse matrix called the {\tmem{reciprocal p-ad}}. The $\tmop{basis} p
- \tmop{ad}$ and its recpiprocal satisify the orthogonality and completeness
relations
\begin{equation}
  e^i_{\mu} e_i^{\nu} = \delta_{\mu}^{\nu}
\end{equation}
and
\begin{equation}
  e_{\mu}^i e_j^{\mu} = \delta_j^i
\end{equation}
The inverse of $\partial_{\mu}$ is
\begin{equation}
  \partial_i = e_i^{\mu} (q) \partial_{\mu}
\end{equation}
which is related to the curvilinear transform of the Cartesian
quantum-mechanical momentum operators by
\begin{equation}
  \hat{p}_i = - i \hbar \partial_i = - i \hbar e_i^{\mu} (q) \partial_{\mu}
\end{equation}
The Hamiltonian operator for free particles is defined by
\begin{equation}
  \hat{H}_0 = \hat{T} = \frac{1}{2 M} \hat{p}^2 = - \frac{\hbar^2}{2 M} \Delta
\end{equation}
where its {\tmem{metric tensor}} is given by
\begin{equation}
  g^{\mu \nu}_{} (q) \equiv e^{\mu}_i (q) e^{\nu}_i (q)
\end{equation}
and its inverse\quad by
\begin{equation}
  g_{\mu \nu} (q) \equiv e^i_{\mu} (q) e^i_{\nu} (q)
\end{equation}
The Laplacian of a metric tensor is then expressed
\begin{equation}
  \begin{array}{ll}
    \Delta & = \partial_i^2\\
    & = e^{i \mu} \partial_{\mu} e^{i \nu} \partial_{\nu}\\
    & = e^{i \mu} e^{i \nu} \partial_{\mu} \partial_{\nu} + (e^{i \mu}
    \partial_{\mu} e^{i \nu}) \partial_{\nu}\\
    & = g^{\mu \nu} (q) \partial_{\mu} \partial_{\nu} - \Gamma_{\mu}^{\mu
    \nu} (q) \partial_{\nu}
  \end{array}
\end{equation}
where
\begin{equation}
  \Gamma_{\mu \nu}^{\lambda} (q) = - e^i_v (q) \partial_{\mu} e_i^{\lambda}
  (q) = e_i^{\lambda} (q) \partial_u e^i_{\nu} (q)
\end{equation}
is the {\tmem{affine-connection}}
{\cite[1.13]{kleinert2004path}}{\cite{podolsky1928}}

\subsubsection{Noether's Theorem and Lie Groups}

From Noether's theorem it is known that continuous symmetries have
corresponding conservation laws.{\cite[4.2.1]{MonsterMoonshine}} Let
\begin{equation}
  \alpha_s
\end{equation}
be a continuous family of symmetries which is a $1$-parameter subgroup $s
\mapsto \alpha_s$ in the Lie group of symmetries. A Lie group is a group whose
operations are compatible with the smooth structure. A smooth structure on a
manifold allows for an unambiguous notion of smooth function.

\subsection{Schroedinger's Time-Dependent Equation and Nonstationary Wave
Motion}

The function $\cos^2 (t) = 1 - \sin^2 (t)$ arises in the Dirac phase averaving
method of calculating transition probabilities of non-stationary states in the
time-dependent Schroedinger equation. {\cite[15.1]{golden1961}} 3

\subsubsection{Operators and Observables: Dirac's Time-Dependent Theory}

Let
\[ y (x, t) = A \cos \left( \frac{2 \pi p x}{h} \right) \cos \left( \frac{2
   \pi \epsilon t}{h} \right) \]
where
\begin{equation}
  \epsilon = h v
\end{equation}
which relates the energies of a ``matter wave'' system, whatever that is,
maybe a fermionic system, to the frequencies of quanta it emits or absorbs.
There is always some arbitriness with choice of units. This wave function has
the symmetry
\begin{equation}
  \frac{\partial^2 y (x, t)}{x^2} = - \left( \frac{2 \pi}{h} \right)^2 p^2 y
  (x, t)
\end{equation}
so that the {\tmem{kinetic energy}} of a particle is obtained from the wave
function
\begin{equation}
  - \frac{(2 \pi)^2}{2 m} \frac{\partial^2 y (x, t)}{x^2} = \frac{p^2}{2 m} y
  (x, t)
\end{equation}
which means formally that
\begin{equation}
  \frac{1}{2 m} (p)^2 y (x, t) = - \frac{1}{2 m} \left( \frac{2 \pi}{i}
  \frac{\partial}{\partial x} \right)^2
\end{equation}
which suggests the identifiation of the {\tmem{Schroedinger momentum
operator}}
\begin{equation}
  \left( \frac{2 \pi}{i} \frac{\partial}{\partial x} \right) \leftrightarrow p
\end{equation}
The equation for ``nonfree'' particles is augmented by
\begin{equation}
  \frac{m}{2} \left( \frac{2 \pi}{i} \frac{\partial}{\partial x} \right)^2 y
  (x, t) = [\epsilon - V (x)] y (x, t)
\end{equation}
where $\epsilon$ is the total energy of the particle and $V (x)$ is its
potential energy. The first Schroedinger equation is then expressed
\begin{equation}
  \epsilon y (x, t) = \left[ \frac{1}{2 m} \left( \frac{2 \pi}{i}
  \frac{\partial}{\partial x} \right)^2 + V (x) \right] y (x, t) = \left[
  \frac{1}{2 m} (p)^2 + V (x) \right] y (x, t)
\end{equation}
or
\[ H_{y (x, t)} = \epsilon y (x, t) \]
where $H$ is the operator corresponding to the Hamiltonian of the (point)
particle. {\cite[11.5]{golden1961}} The {\tmem{time-dependent Schroedinger
equation}} is written
\begin{equation}
  \frac{\partial^2 (x, t)}{\partial t^2} = - \frac{1}{(2 \pi)^2} \epsilon^2 y
  (x, t) = - \left( \frac{H}{2 \pi} \right)^2 y (x, t)
\end{equation}
or equivalently as
\begin{equation}
  0 = \left( (2 \pi)^2 \frac{\partial^2}{\partial t^2} + H^2 \right) y (x, t)
  = \left( H + i 2 \pi \frac{\partial}{\partial t} \right) \left( H - i 2 \pi
  \frac{\partial}{\partial t} \right) y (x, t)
\end{equation}
where
\begin{equation}
  H y (x, t) = i 2 \pi \frac{\partial}{\partial t} y (x, t)
\end{equation}
It is also worth mentioning that in Dirac's theory of the time-dependent
Schroedinger there is a Hamiltonian of the form
\begin{equation}
  H = H_0 + \lambda H_1 (t)
\end{equation}
which has a transition probability per unit time of
\[ w (p \rightarrow - p) = \frac{V_0^2 \sin^2 \left( 2 l \sqrt{\frac{2 m
   \epsilon}{^{\hbar^2}}} \right)}{2 \sqrt{8 m} \epsilon^{\frac{3}{2}} L} \]
where $V_0, \epsilon, h, \tmop{and} L$ are suitable constants.
{\cite[15.5]{golden1961}}

\subsubsection{The String Theoretic Partition Function}

The {\tmem{string theoretic partition function}} is defined as
\begin{equation}
  Z_{s t} (\omega) = Z_{s t, R} (\omega) = \tmop{Tr} [e^{2 \pi i \omega_1 P -
  2 \pi \omega_2 H}] = (q \bar{q})^{- \frac{1}{24}} \tmop{Tr} (q^{L_0^+}
  \bar{q}^{L^-_0})
\end{equation}
where $L^+_0$ and $L^-_0$ are the Virasoro generators defined by
\begin{equation}
  L_0^{\pm} = \frac{(p^+)^2 + (p^-)^2}{2} + \sum_{n = 1}^{\infty} \alpha_{-
  n}^{\pm} \alpha_n^{\pm}
\end{equation}
and $\alpha_{- n}^{\pm} \alpha_n$ are related to something called a
Fubini-Veneziano field and $p^{_+}$ is the {\tmem{left}}-momentum operator and
$p^-$ is {\tmem{right-momentum}} operator. {\cite[2.2.1]{fractalzetastrings}}

\


\begin{thebibliography}{Coo03}
  \bibitem[11]{prz}D.~Schumayer  and  D.~A.~W.~Hutchinson.{\newblock}
  Colloquium: Physics of the Riemann hypothesis.{\newblock} \tmtextit{Reviews
  of Modern Physics}, 83:307--330, apr 2011.{\newblock}
  
  \bibitem[CB90]{cva}Ruel~V.~Churchill  and  James~Ward Brown.{\newblock}
  \tmtextit{Complex Variables and Applications}.{\newblock} McGraw-Hill
  Publishing Company, 1990.{\newblock}
  
  \bibitem[Coo03]{Psm}David~B.~Cook.{\newblock} \tmtextit{Probability and
  Schrodinger's Mechanics}.{\newblock} World Scientific Publishing Company,
  2003.{\newblock}
  
  \bibitem[Fl{\"u}94]{pqm}S.~Fl{\"u}gge.{\newblock} \tmtextit{Practical
  Quantum Mechanics}.{\newblock} Classics in Mathematics. Springer Berlin
  Heidelberg, 1994.{\newblock}
  
  \bibitem[Gan06]{MonsterMoonshine}Terry Gannon.{\newblock}
  \tmtextit{Moonshine Beyond the Monster: The Bridge Connecting Algebra,
  Modular Forms, and Physics}.{\newblock} Cambridge University Press,
  Monographs on Mathematical Physics, 2006.{\newblock}
  
  \bibitem[Gol61]{golden1961}Sydney Golden.{\newblock} \tmtextit{An
  Introduction to Theoretical Physical Chemistry}.{\newblock} Series In
  Chemistry. Addison-Wesley Publishing Company, Inc., 1961.{\newblock}
  
  \bibitem[Ivi13]{HardyZ}A.~Ivi{\'c}.{\newblock} \tmtextit{The Theory of
  Hardy's Z-Function}.{\newblock} Cambridge Tracts in Mathematics. Cambridge
  University Press, 2013.{\newblock}
  
  \bibitem[Kle04]{kleinert2004path}H.~Kleinert.{\newblock} \tmtextit{Path
  Integrals in Quantum Mechanics, Statistics, Polymer Physics, and Financial
  Markets}.{\newblock} World Scientific, 2004.{\newblock}
  
  \bibitem[Kna99]{knauf99}Andreas Knauf.{\newblock} Number theory, dynamical
  systems and statistical mechanics.{\newblock} \tmtextit{Reviews in
  Mathematical Physics}, 11(08):1027--1060, 1999.{\newblock}
  
  \bibitem[Lap08]{fractalzetastrings}Michel~L.~Lapidus.{\newblock}
  \tmtextit{In search of the Riemann zeros: Strings, Fractal membranes and
  Noncommutative Spacetimes}.{\newblock} American Mathematical Society,
  2008.{\newblock}
  
  \bibitem[Pod28]{podolsky1928}Boris Podolsky.{\newblock} Quantum-mechanically
  correct form of hamiltonian function for conservative systems.{\newblock}
  \tmtextit{Phys. Rev.}, 32:812--816, Nov 1928.{\newblock}
\end{thebibliography}
\end{document}